\documentclass[a4paper,12pt]{amsart}
\usepackage{amssymb}
\usepackage{ifthen}
\usepackage[usenames]{color}
\usepackage{graphicx}
\nonstopmode \numberwithin{equation}{section}
\newtheorem{thm}{Theorem}[section]

\newtheorem{cor}{Corollary}[section]
\newtheorem{lem}{Lemma}[section]
\newtheorem{prop}{Proposition}[section]

\newtheorem{conj}[equation]{Conjecture}

\theoremstyle{definition}
\newtheorem{defn}{Definition}[section]
\newtheorem{examp}{Example}[section]
\newtheorem{prob}[equation]{Problem}
\newtheorem{ques}[equation]{Question}
\newtheorem{rem}{Remark}[section]

\newcounter {own}
\def\theown {\thesection       .\arabic{own}}

\newenvironment{pf}[1][]{%
 \vskip 3mm
 \noindent
 \ifthenelse{\equal{#1}{}}%
  {{\slshape Proof. }}%
  {{\slshape #1.} }%
 }%
{\qed\bigskip}

\newcounter{alphabet}
\newcounter{tmp}
\newenvironment{Thm}[1][]{\refstepcounter{alphabet}%
\bigskip%
\noindent%
{\bf Theorem \Alph{alphabet}}%
\ifthenelse{\equal{#1}{}}{}{ (#1)}%
{\bf .} \itshape}{\vskip 8pt}

\makeatletter
\newcommand{\Ref}[1]{\@ifundefined{r@#1}{}{\setcounter{tmp}{\ref{#1}}\Alph{tmp}}}
\makeatother

\newenvironment{Lem}[1][]{\refstepcounter{alphabet}%
\bigskip%
\noindent%
{\bf Lemma \Alph{alphabet}}%
{\bf .} \itshape}{\vskip 8pt}

\newcommand{\IR}{{\mathbb R}}

\newcommand{\IC}{{\mathbb C}}

\newcommand{\ID}{{\mathbb D}}

\def\be{\begin{equation}}
\def\ee{\end{equation}}

\newcommand{\bee}{\begin{enumerate}}
\newcommand{\eee}{\end{enumerate}}

\newcommand{\blem}{\begin{lem}}
\newcommand{\elem}{\end{lem}}
\newcommand{\bthm}{\begin{thm}}
\newcommand{\ethm}{\end{thm}}
\newcommand{\bcor}{\begin{cor}}
\newcommand{\ecor}{\end{cor}}
\newcommand{\beg}{\begin{examp}}
\newcommand{\eeg}{\end{examp}}
\newcommand{\begs}{\begin{examples}}
\newcommand{\eegs}{\end{examples}}
\newcommand{\bdefe}{\begin{defn}}
\newcommand{\edefe}{\end{defn}}
\newcommand{\bprob}{\begin{prob}}
\newcommand{\eprob}{\end{prob}}
\newcommand{\bques}{\begin{ques}}
\newcommand{\eques}{\end{ques}}
\newcommand{\bei}{\begin{itemize}}
\newcommand{\eei}{\end{itemize}}

\newcommand{\bcon}{\begin{conj}}
\newcommand{\econ}{\end{conj}}
\newcommand{\bcons}{\begin{conjs}}
\newcommand{\econs}{\end{conjs}}
\newcommand{\bprop}{\begin{prop}}
\newcommand{\eprop}{\end{prop}}
\newcommand{\br}{\begin{rem}}
\newcommand{\er}{\end{rem}}
\newcommand{\brs}{\begin{rems}}
\newcommand{\ers}{\end{rems}}
\newcommand{\bo}{\begin{obser}}
\newcommand{\eo}{\end{obser}}
\newcommand{\bos}{\begin{obsers}}
\newcommand{\eos}{\end{obsers}}
\newcommand{\bpf}{\begin{pf}}
\newcommand{\epf}{\end{pf}}
\newcommand{\ba}{\begin{array}}
\newcommand{\ea}{\end{array}}
\newcommand{\beq}{\begin{eqnarray}}
\newcommand{\beqq}{\begin{eqnarray*}}
\newcommand{\eeq}{\end{eqnarray}}
\newcommand{\eeqq}{\end{eqnarray*}}

\newcounter{minutes}
\divide\time by 60
\newcounter{hours}
\multiply\time by 60 \addtocounter{minutes}{-\time}

\begin{document}
\bibliographystyle{amsplain}
\title [Lipschitz continuity of quasiconformal mappings]
{Lipschitz continuity of quasiconformal mappings and of the solutions to second order elliptic PDE
 with respect to the distance ratio metric}

\def\thefootnote{}
\footnotetext{ \texttt{\tiny File:~\jobname .tex,
          printed: \number\day-\number\month-\number\year,
          \thehours.\ifnum\theminutes<10{0}\fi\theminutes}
} \makeatletter\def\thefootnote{\@arabic\c@footnote}\makeatother

\author{Peijin Li}
\address{Peijin Li, Department of Mathematics,
Hunan First Normal University, Changsha, Hunan 410205,
People's Republic of China}
\email{wokeyi99@163.com}
\author{Saminathan Ponnusamy}
\address{S. Ponnusamy,  Stat-Math Unit,\\
Indian Statistical Institute (ISI), Chennai Centre,
110, Nelson Manickam Road, Aminjikarai,\\
Chennai, 600 029, India. }
\email{samy@isichennai.res.in, samy@iitm.ac.in}

\subjclass[2010]{Primary: 30C62, 30L10, 31A30,   31B30; Secondary: 26A16, 33C05, 35J05, 30C20}
\keywords{Quasiconformal mappings, harmonic and polyharmonic mappings, distance ratio metric, and Lipschitz continuity.
}

\begin{abstract}
The main aim of this paper is to study the Lipschitz continuity of certain $(K, K')$-quasiconformal mappings
with respect to the distance ratio metric, and the Lipschitz continuity of the solution of a quasilinear differential
equation with respect to the distance ratio metric.
\end{abstract}


\maketitle \pagestyle{myheadings} \markboth{P. Li and S. Ponnusamy
}{Lipschitz continuity of quasiconformal mappings}

\section{Introduction and main results}\label{csw-sec1}

Martio \cite{M} was the first who considered the study on harmonic quasiconformal mappings in $\IC$. In the recent years
the articles  \cite{chen, K1,KM0, KM1, KP0, P} brought much light
on this topic. In \cite{CLSW, LCW}, the Lipschitz characteristic of $(K, K')$-quasiconformal mappings
has been discussed. In \cite{MKM}, the authors proved that a $K$-quasiconformal harmonic mapping from the unit disk $\ID$
onto itself is bi-Lipschitz with respect to hyperbolic metric, and also proved that a $K$-quasiconformal harmonic mapping
from the upper half-plane $\mathbb{H}$ onto itself is bi-Lipschitz with respect to hyperbolic metric. In \cite{MV},
the authors proved that a $K$-quasiconformal harmonic mapping from $D$ to $D'$
is bi-Lipschitz with respect to quasihyperbolic metrics on $D$ and $D'$, where $D$ and $D'$ are proper domains in $\IC$.
Important definitions will be included later in this section.


In \cite{K2}, Kalaj considered the bi-Lipschitz continuity of $K$-quasiconformal solution of the inequality
\be\label{eq-2.1}
|\Delta f|\leq B|D f|^2.
\ee
Here $\Delta f$ represents the two-dimensional Laplacian of $f$ defined by
$\Delta f=f_{xx}+f_{yy}=4f_{z\overline{z}}
$
and the mapping $f$ satisfying the Laplace equation $\Delta f=0$ is called harmonic.
For $z=x+iy$ and $f=u+iv,$ $Df$ denotes the Jacobian matrix
$$\left(
  \begin{array}{cc}
    u_x & u_y \\
    v_x & v_y \\
  \end{array}
\right)
$$
so that  $J_f=|f_z|^2-|f_{\overline z}|^2$ is the Jacobian of $f$.

The first aim of this paper is to consider the Lipschitz continuity of $(K, K')$-quasiconformal solution of
the inequality \eqref{eq-2.1} with respect to the distance ratio metric.

\begin{thm}\label{thm2.4}
Let $f$ be a $(K, K')$-quasiconformal $C^2$ mapping from the unit disk $\mathbb{D}=\{z:\,|z|<1\}$ onto itself,
satisfying the inequality \eqref{eq-2.1} and $f(0)=0$. Then $f$ is Lipschitz continuous with respect to the distance ratio metric.
\end{thm}

The proof of Theorem \ref{thm2.4} will be presented in Section \ref{sec-2}.
Before we proceed further, let us fix up further notation, preliminaries and remarks. The proof of Theorem \ref{thm2.4} will be presented in
Section \ref{sec-2}.

\subsection{$(K, K')$-quasiconformal mappings}

We say that a function $u:\, D\rightarrow \mathbb{R}$ is {\it absolutely continuous on lines}, {\it ACL} in brief,
in the domain $D$ if for every closed rectangle $R\subset D$ with sides parallel to the axes $x$ and $y$,
$u$ is absolutely continuous on almost every horizontal line segment and almost every vertical line segment in $R$.
Such a function has, of course, partial derivatives $u_x$ and $u_y$ a.e. in $D$ (cf. \cite{Ah}).
Further, we say $u\in ACL^{2}$ if $u\in ACL$ and its partial derivatives are locally $L^2$ integrable in $D$.

A sense-preserving continuous mapping $f:\, D\rightarrow \Omega $ is said to be
\begin{enumerate}
\item {\it $(K, K')$-quasiregular} if $f$ is $ACL^{2}$ in $D$, $J_{f}\not=0$
a.e. in $D$ and there are constants $K\geq1$ and $K'\geq0$
such that $|D f|^2\leq KJ_f+K'$ a.e. in $D$, where $|D f|=|f_z|+|f_{\overline z}|$;
\item {\it $K$-quasiregular} if $K'=0$.
\end{enumerate}

In particular, $f$ is called {\it $(K, K')$-quasiconformal} if $f$ is a $(K, K')$-quasiregular homeomorphism;
and $f$ is {\it $K$-quasiconformal} if $f$ is a $K$-quasiregular homeomorphism.

Here are some basic comments on these mappings.
From \cite[Example 2.1]{CLSW} and \cite[Example 2.1]{LCW} we know that there are $(K, K')$-quasiregular mappings which are not
$K_1$-quasiregular for any $K_1\geq 1$. Moreover, it is known that (see \cite[Example 4.1]{CLSW}) there are
$(K, K')$-quasiconformal mappings whose inverses are not $(K_1, K_1')$-quasiconformal for any $K_1\geq 1$ and $K_1'\geq 0$.

\br
If $f$ is a $(K, K')$-quasiregular mapping, $g$ is a analytic function and $|g'|$ is bounded by a constant $L$,
then $f\circ g$ is $(K, K'L^2)$-quasiregular mapping.
\er



A mapping $f:\,D\to \Omega$ is {\it proper} if the preimage of every compact set in $\Omega$ is compact in $D$
 (cf. \cite[p.~4051]{KP} or \cite[p.~17]{V}).

\subsection{The distance ratio metric}

For a subdomain $G\subset \mathbb{C}$ and for all $z$, $w\in G$, the distance ratio metric $j_{G}$ is defined as
$$j_{G}(z,w)=\log\left(1+\frac{|z-w|}{\min\{\delta_G(z),\delta_G(w)\}}\right),$$
where $\delta_G(z)$ denotes the Euclidean distance from $z$ to $\partial G$. The distance ratio metric
was introduced by Gehring and Palka \cite{go2} and in the above simplified form by Vuorinen \cite{vo2}.
However, the distance ratio metric $j_{G}$ is not invariant under M\"{o}bius transformations. Therefore, it
is natural to consider the Lipschitz continuity of  conformal mappings or M\"{o}bius transformations with respect
to the distance ratio metric. Gehring and Osgood \cite{go1} proved that the distance ratio metric is not altered
by more than a factor of $2$ under M\"{o}bius transformations. 

\begin{Thm}\label{ThmA}{\rm (\cite[Proof of Theorem 4]{go1})}
 If $G$ and $G'$ are proper subdomains of $\mathbb{R}^{n}$ and if $f$ is a M\"{o}bius transformation of $G$ onto $G'$, then
$j_{G'}(f(x),f(y))\leq 2j_{G}(x,y)$  for all $x$, $y\in G$.
\end{Thm}


Recall that a mapping $f:\,D\to \Omega$ is said to be {\it Lipschitz continuous}  (resp. Lipschitz continuous with respect to the distance ratio metric)
if there exists a positive constant $L_{1}$ (resp. a positive constant $L$) such that for all $z, w\in D$,
$$|f(z)-f(w)|\leq L_{1}|z-w| \quad \mbox{(resp. $j_{\Omega}(f(z),f(w))\leq L j_{D}(z,w)$)}.
$$


In 2011, Kalaj and  Mateljevi\'{c} \cite{KM1} proved that every quasiconformal $C^2$ diffeomorphism $f$
from the domain $\Omega$ with $C^{1, \alpha}$ compact boundary onto the domain $G$ with $C^{2, \alpha}$ compact boundary
satisfying the {\it Poisson differential inequality}
\be\label{eq-1}
|\Delta f|\leq B|D f|^2+C
\ee
for some constants $B\geq 0$ and $C\geq 0$, is Lipschitz continuous respect to Euclidean metric.
Clearly if $B=C=0$ in \eqref{eq-1}, then $f$ is harmonic.

Recently, the authors in \cite[Theorem 1.1]{CLSW} proved the following theorem:

\begin{Thm}\label{ThmB}
Suppose $f$ is a proper $(K, K')$-quasiregular $C^2$ mapping of a Jordan domain $D$
with $C^{1,\alpha}$ boundary onto a Jordan domain $\Omega$ with $C^{2,\alpha} $ boundary. If $f$ satisfies the partial differential inequality \eqref{eq-1} for constants $B>0$ and $C\geq 0$,
then $f$ has bounded partial derivatives in $D$.
In particular, $f$ is Lipschitz continuous.
\end{Thm}

\br
From Theorem \Ref{ThmB} we  infer that if $f:\,\ID \to \ID$ satisfies the conditions of Theorem~\ref{thm2.4}, then there exists a constant
$M$ such that $|D f|\leq M$. Hence for all $z, w\in \ID$, we have $|f(z)-f(w)|\leq M|z-w|$ and $|f(z)|\leq M|z|$.
If $M <1$,  we get
 \begin{align*}
j_{\ID}(f(z), f(w))=&\log\left(1+\frac{|f(z)-f(w)|}{\min\{\delta_{\ID}(f(z)),\delta_{\ID}(f(w))\}}\right)\\
\leq&\log\left(1+\frac{M|z-w|}{M\,\min\{\delta_{\ID}(z),\delta_{\ID}(w)\}}\right)\\
\leq&j_{\ID}(z, w),
\end{align*}
which clearly shows that $f:\,\ID \to \ID$ is a Lipschitz continuous function with respect to the distance ratio metric. \hfill $\Box$
\er


%

In order to state our next result,
we need to recall the definition of hypergeometric series.
For $a, b, c\in\IR$ with $c\neq 0, -1, -2, \ldots$, the {\it hypergeometric} function is defined by the power series
$$
F(a, b; c; z):={}_2F_1(a,b;c;z)=\sum_{n=0}^{\infty}\frac{(a)_n(b)_n}{(c)_n}\frac{z^n}{n!},\;\;|z|<1,
$$
where $(a)_0=1$ and $(a)_n=a(a+1)\cdots(a+n-1)$ for $n=1, 2, \ldots$ are the {\it Pochhammer symbols}. Obviously, for $n=0, 1, 2, \ldots$,
$(a)_n=\Gamma(a+n)/\Gamma(a)$. In particular, for $a, b, c>0$ and $a+b<c$, we have (cf. \cite{AB, AQ})
$$F(a, b; c; 1)=\lim_{z\rightarrow 1^{-}}F(a, b; c; z)=\frac{\Gamma(c)\Gamma(c-a-b)}{\Gamma(c-a)\Gamma(c-b)}<\infty.
$$

Consider the operator equation
\be\label{eq-11}
 T_{\alpha}(f)=0 ~ \mbox{ in $\ID$},
\ee
where $f:\,\ID \to \IC$, $\alpha\in\IR$, and
$$T_\alpha=-\frac{\alpha^2}{4}(1-|z|^2)^{-\alpha-1}
+\frac{\alpha}{2}(1-|z|^2)^{-\alpha-1}(z\frac{\partial}{\partial z}+\overline{z}\frac{\partial}{\partial \overline{z}})
+(1-|z|^2)^{-\alpha}\frac{\partial^2}{\partial z\partial \overline{z}}
$$
is the {\it second order elliptic partial differential operator} defined on the unit disk $\ID$.  In the case of $\alpha=0$,
$T_{\alpha}(f)=0$ is equivalent to saying that $f$ is harmonic in $\ID$. More generally, if $f$ satisfies \eqref{eq-11} with $\alpha=2(n-1)$,
then $f$ is {\it polyharmonic} (or {\it $n$-harmonic}) in $\ID$, where $n\in\{1, 2, \ldots\}$ (cf. \cite{AGPon-17, BH, CPW, P1}).
Recently, several new properties of polyharmonic mappings are discussed in \cite{AGPon-17}.
The following result concerns the solutions to the equation \eqref{eq-11}.

\begin{Lem}\label{LemF}\cite[Theorem 2.2]{O}
Let $\alpha\in\IR$ and $f\in C^2(\ID)$. Then $f$ satisfies \eqref{eq-11} if and only if it has a series expansion for $z\in \ID$ of the form
\be
\label{eq-31}
f(z)=\sum_{k=0}^{\infty}c_kF(-\frac{\alpha}{2}, k-\frac{\alpha}{2}; k+1; |z|^2)z^k
+\sum_{k=1}^{\infty}c_{-k}F(-\frac{\alpha}{2}, k-\frac{\alpha}{2}; k+1; |z|^2){\overline{z}}^k,
\ee
%
where $\{c_k\}_{k=-\infty}^\infty$ is a sequence of complex numbers satisfying
\be\label{eq-32}
\limsup_{|k|\rightarrow\infty}|c_k|^{\frac{1}{|k|}}\leq1.
\ee
In particular, the expansion \eqref{eq-31}, subject to \eqref{eq-32}, converges in $C^\infty(\ID)$, and every solution $f$ of \eqref{eq-11}
is $C^\infty$-smooth in $\ID$.
\end{Lem}

In \cite{CV, MC}, the authors gave some properties of solution to (\ref{eq-11}) whereas in \cite{si, sivw}, the authors considered the
Lipschitz continuity of the distance-ratio metric under some M\"{o}bius automorphisms of the unit ball and conformal mappings from $\mathbb{D}$ to $\mathbb{D}$.
In \cite{CRW}, the authors discussed the Lipschitz continuity of polyharmonic mappings with respect to the distance ratio metric.
Thus, it is natural to investigate Lipschitz continuity of the solution of (\ref{eq-11}) in $\ID$ with respect to the
distance ratio metric. We now state our next result.

\begin{thm}\label{thm3.1}
For $\alpha>-1$, let $f:\, \ID \rightarrow \ID$ be a  $C^2$-solution to \eqref{eq-11} with the series expansion
of the form \eqref{eq-31} and $f(0)=0$. If
\be\label{eq31}
\sum_{k=1}^{\infty}(|c_{k}|+|c_{-k}|)\sum_{n=0}^{\infty}\frac{(-\frac{\alpha}{2})_n(k-\frac{\alpha}{2})_n}{(k+1)_n n!}\leq 1,
\ee
then $$j_{\mathbb{D}}(f(z),f(w))\leq j_{\mathbb{D}}(z,w),$$
and this inequality is sharp.
That is, $f$ is Lipschitz continuous with respect to the distance ratio metric.
\end{thm}

\br
In Theorem \ref{thm3.1}, we restrict $\alpha>-1$, see \cite[proposition 1.4]{O} for the reason for this constraint.
\er

The proof of Theorem \ref{thm3.1} will be presented in Section \ref{sec-2}.

\section{Proof of Theorems \ref{thm2.4} and  \ref{thm3.1}}\label{sec-2}

First we shall deal with the Lipschitz continuity of certain $(K, K')$-quasiconformal mappings and then
consider the Lipschitz continuity of the solution  to the differential operator $T_\alpha$ with respect to the
distance ratio metric.

\begin{lem}\label{lem2}
Assume the hypotheses of Theorem {\rm \ref{thm2.4}}. Then there exists a constant $C(K, K', B)$ such that for $z\in\ID$
\be\label{eq-2.2}
\frac{1-|z|^2}{1-|f(z)|^2}\leq C(B,K, K').
\ee
\end{lem}
\bpf By assumption $f$ is a $(K, K')$-quasiconformal $C^2$ mapping from $\ID$ onto itself,
satisfying the inequality \eqref{eq-2.1} and $f(0)=0$. For convenience, we denote the class of all such
functions $f$ by $\mathcal{QC}(\ID, B, K, K')$.
Then there is a positive constant $A$ not depending on $f$ such that the function $\varphi_f$, $f\in \mathcal{QC}(\ID, B, K, K')$,
defined by
$$\varphi_f(z)=-\frac{1}{A}+\frac{1}{A}e^{A(|f(z)|-1)}
$$
is subharmonic in $\ID$.

Now, let us prove the existence of such an $A$. Take
$$\psi(\rho)=-\frac{1}{A}+\frac{1}{A}e^{A(\rho-1)}.$$
Then $\psi'(\rho)=e^{A(\rho-1)}$ and $\psi''(\rho)=Ae^{A(\rho-1)}$. On the other hand, using
$f_z=(1/2)(f_x-if_y)$ and $f_{\overline{z}} =(1/2)(f_x+if_y)$, we find that
$$|D |f||^2=|f|_x^2+|f|_y^2=(|f|_z+|f|_{\overline{z}})^2+i^2(|f|_z-|f|_{\overline{z}})^2=4|f|_z|f|_{\overline{z}},
$$
and thus,
\be\label{eq-2.3}
\Delta\varphi_f=\psi''(|f|)|D |f||^2+\psi'(|f|)\Delta |f|.
\ee
Furthermore, put $s=f/|f|$. By elementary calculations we see that the following equalities hold:
$$\Delta|f|=\frac{|f_z|^2+|f_{\overline{z}}|^2}{|f|}-2{\rm Re}\left(\overline{f}^{\frac{1}{2}}f^{-\frac{3}{2}}f_zf_{\overline{z}}\right)+{\rm Re}\,(\overline{s}\Delta f)$$
and
$$|s_z|^2=|s_{\overline{z}}|^2
=\frac{1}{4}\frac{|f_z|^2+|f_{\overline{z}}|^2}{|f|^2}-\frac{1}{2|f|}{\rm Re}\left(\overline{f}^{\frac{1}{2}}f^{-\frac{3}{2}}f_zf_{\overline{z}}\right).$$
Then we know that
\be\label{eq-2.4}
\Delta |f|=|f|\cdot|D s|^2+{\rm Re}\,(\overline{s}\Delta f).
\ee
We continue the discussion by setting $\rho=|f|$. According to \cite[Lemma 3.1]{CLSW}, we have
\be\label{eq-2.5}
|D \rho|\geq\frac{|D f|}{K}-\frac{\sqrt{K'}}{K},
\ee
Using \eqref{eq-2.1}, \eqref{eq-2.3}, \eqref{eq-2.4} and \eqref{eq-2.5}, it follows finally that
\begin{align*}
\Delta\varphi_f=&e^{A(\rho-1)}\Big[A|D \rho|^2+\rho|D s|^2+{\rm Re}\,(\overline{s}\Delta f)\Big]\\
\geq& e^{A(\rho-1)}\Big[\frac{A}{K^2}(|D f|-\sqrt{K'})^2-B|D f|^2\Big]\\
=&e^{A(\rho-1)}\Big(\frac{A-BK^2}{K^2}|D f|^2-\frac{2A\sqrt{K'}}{K^2}|D f|+\frac{AK'}{K^2}\Big).
\end{align*}
We obtain from Theorem \Ref{ThmB} that $f$ is Lipschitz continuous, and then there exists a constant $M$ such
that $|D f|\leq M$. Hence, if we choose an appropriate $A$, satisfying $A\neq BK^2$ and
$$\frac{A\sqrt{K'}-\sqrt{ABK'K^2}}{A-BK^2}\geq M,
$$
i.e.
\be\label{eq-00}
(\sqrt{K'}-M)^2A^2+[2BMK^2(\sqrt{K'}-M)-BK^2K']A+B^2M^2K^4\geq0,
\ee
we obtain the inequality $\Delta\varphi_f(z)\geq 0$  for $|z|<1$.

If $M=\sqrt{K'}$, then there exists an appropriate value of $A$ satisfying the inequality \eqref{eq-00}.
If $M\neq\sqrt{K'}$ and $K'+4M^2-4M\sqrt{K'}\leq0$, then \eqref{eq-00} holds for all $A$.
If $M\neq\sqrt{K'}$ and $K'+4M^2-4M\sqrt{K'}>0$, then \eqref{eq-00} holds for all
$$A\geq\frac{BK^2K'-2BMK^2(\sqrt{K'}-M)+BK^2\sqrt{K'(K'+4M^2-4M\sqrt{K'})}}{2(\sqrt{K'}-M)^2}.$$
In conclusion, there must exist an appropriate $A$ such that $\Delta\varphi_f(z)\geq 0$  for $|z|<1$.

Define
$$F(z)=\sup\{\varphi_f(z):\, f\in\mathcal{QC}(\ID, B, K, K')\}.$$
We prove that $F$ is subharmonic in $\ID$. By \cite[Theorem 1.6.2]{H1}, we only need to prove that $F$ is continuous.
Define
$h(z)=e^{A(|z|-1)}, \;|z|<1.$
Elementary calculations show that
$$h_z(z)=\frac{A}{2}\frac{\overline{z}}{|z|}e^{A(|z|-1)}
\;\;\mbox{and}\;\;h_{\overline{z}}(z)=\frac{A}{2}\frac{z}{|z|}e^{A(|z|-1)}.$$
Then $|D h|=|h_z|+|h_{\overline z}|=Ae^{A(|z|-1)}<A$ which implies that
$$|h(z)-h(z')|\leq A|z-z'| ~\mbox{ for $z, z'\in\ID$}.
$$
According to Theorem \Ref{ThmB}, we know that $f$ is Lipschitz continuous. Therefore
$$|\varphi_f(z)-\varphi_f(z')|= \frac{1}{A}\left |e^{A(|f(z)|-1)}-e^{A(|f(z')|-1)}\right |\leq|f(z)-f(z')| \leq M|z-z'|,
$$
where $M$ is a constant. Hence, $|F(z)-F(z')|\leq M|z-z'|$ so that $F$ is continuous. Finally,
from the similar proof of \cite[Lemma 2.3]{K2}, we complete the proof.
\epf

\subsection{ Proof of Theorem \ref{thm2.4}}
From the hypotheses of Theorem \ref{thm2.4} and Lemma~\ref{lem2}, we obtain that
$$\frac{1-|z|^2}{1-|f(z)|^2}\leq C(K, K', B)
$$
and thus, we obtain that
$$\frac{1-|z|}{1-|f(z)|}\leq C(K, K', B)\frac{1+|f(z)|}{1+|z|}\leq2C(K, K', B).
$$
Moreover, from Theorem \Ref{ThmB}, we see that $f$ is Lipschitz continuous and therefore, there exists a constant
$M_1$ such that $|D f|\leq M_1$. Now, we choose an appropriate constant $M$ satisfying $M>\max\{M_1, 1/(2C(K, K', B))\}$
so that $|D f|\leq M$. Consequently, using the Bernoulli inequality,
for any two points $z$ and $w$ in $\ID$, we have
\begin{align*}
j_{\ID}(f(z), f(w))=&\log\left(1+\frac{|f(z)-f(w)|}{\min\{\delta_{\ID}(f(z)),\delta_{\ID}(f(w))\}}\right)\\
\leq&\log\left(1+2C(K, K', B)M\frac{|z-w|}{\min\{\delta_{\ID}(z),\delta_{\ID}(w)\}}\right)\\
\leq&2C(K, K', B)Mj_{\ID}(z, w)
\end{align*}
and thus, the proof of the theorem is complete. \hfill $\Box$




\subsection{Proof of Theorem \ref{thm3.1}}
For convenience, let
$g(t)=F(-\frac{\alpha}{2}, k-\frac{\alpha}{2}; k+1; t).
$
For $z,w\in \mathbb{D}$, let us assume that $|f(z)|\geq |f(w)|$. Then

\vspace{8pt}
\noindent
$|f(z)-f(w)|$
\begin{align*}
=&\left |\sum_{k=1}^{\infty}c_{k}(|g(|z|^2)z^k-g(|w|^2)w^k)+\sum_{k=1}^{\infty}c_{-k}(|g(|z|^2){\overline{z}}^k-g(|w|^2){\overline{w}}^k)\right |\\
\leq&|z-w|\sum_{k=1}^{\infty}\frac{|g(|z|^2)z^k-g(|z|^2)w^k+g(|z|^2)w^k-g(|w|^2)w^k|}{|z-w|}|c_{k}|\\
&+|z-w|\sum_{k=1}^{\infty}\frac{|g(|z|^2){\overline{z}}^k-g(|z|^2){\overline{w}}^k+g(|z|^2){\overline{w}}^k-g(|w|^2){\overline{w}}^k|}{|z-w|}|c_{-k}|\\
\leq&|z-w|\sum_{k=1}^{\infty}(|c_{k}|+|c_{-k}|)\Big(\sum_{n=0}^{\infty}\frac{(-\frac{\alpha}{2})_n(k-\frac{\alpha}{2})_n}{(k+1)_n n!}|z|^{2n}(|z|^{k-1}+|z|^{k-2}|w|+\cdots+|w|^{k-1})\\
&+\sum_{n=0}^{\infty}\frac{(-\frac{\alpha}{2})_n(k-\frac{\alpha}{2})_n}{(k+1)_n n!}|w|^k(|z|^{2n-1}+|z|^{2n-2}|w|+\cdots+|w|^{2n-1})\Big)\\
\leq&|z-w|\sum_{k=1}^{\infty}(|c_{k}|+|c_{-k}|)\sum_{n=0}^{\infty}\frac{(-\frac{\alpha}{2})_n(k-\frac{\alpha}{2})_n}{(k+1)_n
n!}\sum_{s=0}^{2n+k-1}|z|^s,
\end{align*}
and
\begin{align*}
1-|f(z)|\geq &\sum_{k=1}^{\infty}(|c_{k}|+|c_{-k}|)\sum_{n=0}^{\infty}\frac{(-\frac{\alpha}{2})_n(k-\frac{\alpha}{2})_n}{(k+1)_n n!}-|f(z)|\\
\geq &\sum_{k=1}^{\infty}(|c_{k}|+|c_{-k}|)\Big(\sum_{n=0}^{\infty}\frac{(-\frac{\alpha}{2})_n(k-\frac{\alpha}{2})_n}{(k+1)_n n!}-g(|z|^2)|z|^k\Big)\\
=&\sum_{k=1}^{\infty}(|c_{k}|+|c_{-k}|)\sum_{n=0}^{\infty}\frac{(-\frac{\alpha}{2})_n(k-\frac{\alpha}{2})_n}{(k+1)_n n!}(1-|z|^{2n+k})\\
=&(1-|z|)\sum_{k=1}^{\infty}(|c_{k}|+|c_{-k}|)\sum_{n=0}^{\infty}\frac{(-\frac{\alpha}{2})_n(k-\frac{\alpha}{2})_n}{(k+1)_n n!}\sum_{s=0}^{2n+k-1}|z|^s,
\end{align*}
so that,  using the Bernoulli inequality, we have
\begin{align*}
j_{\mathbb{D}}(f(z),f(w))=&\log\left(1+\frac{|f(z)-f(w)|}{1-|f(z)|}\right)\\
\leq&\log\left(1+\frac{|z-w|}{1-|z|}\right)\\
\leq&j_{\mathbb{D}}(z,w).
\end{align*}
As in \cite[Theorem 7]{CRW}, the mapping $f(z)= |z|^{2(p-1)}z^{m}$ or $f(z)= |z|^{2(p-1)}\overline{z}^{m}$ for $p$, $m\geq1$, shows the sharpness
in the last inequality. The proof of the theorem is complete.  \hfill $\Box$

\subsection*{Acknowledgements}
The authors thank the referee for his/her careful reading and many useful comments.
The first author was supported by Centre for International Co-operation in Science (CICS)
through the award of ``INSA JRD-TATA Fellowship"  and was completed
during her visit to the Indian Statistical Institute (ISI), Chennai Centre.
The research was partly supported by NSF of China (No. 11571216 and No. 11671127).
The second author is on leave from the IIT Madras.



\begin{thebibliography}{99}


\bibitem{Ah} L.V. Ahlfors,
Lectures on quasiconformal mappings, Princeton (NJ), Van Nostrand Company,
1966.

\bibitem{AGPon-17}  K. F. Amozova, E. G. Ganenkova, and S. Ponnusamy,
Criteria of univalence and fully $\alpha$-accessibility for $p$-harmonic and $p$-analytic functions,
\textit{Complex Var. Elliptic Equ.}, \textbf{62}(8) (2017), 1165--1183.

\bibitem{AB} G. Anderson, R. Barnard, K. Richards, M. Vamanamurthy and M. Vuorinen,
Inequalities for zero-balanced hypergeometric functions,
\textit{Trans. Amer. Math. Soc.}, {\bf 126} (1995), 1713--1723.

\bibitem{AQ} G. Anderson, S. Qiu, M. Vamanamurthy and M. Vuorinen,
Generalized elliptic integrals and modular equations,
\textit{Pacific J. Math.}, {\bf 192} (2000), 1--37.





\bibitem{BH} A. Borichev and H. Hedenmalm,
Weighted integrability of polyharmonic functions,
\textit{Adv. Math.}, {\bf 264} (2014), 464--505.



\bibitem{CLSW} J. Chen, P. Li, S. K. Sahoo and X. Wang,
On the Lipschitz continuity of certain quasiregular mappings between smooth Jordan domains,
\textit{Israel J. Math. }, (2017), 26 pages.\\
{\tt DOI 10.1007/s11856-017-1522-y.}

\bibitem{CRW} J. Chen, A. Rasila and X. Wang,
On lengths, areas and Lipschitz continuity of polyharmonic mappings,
\textit{J. Math. Anal. Appl.}, {\bf 422} (2015), 1196--1212.

\bibitem{chen} M. Chen and X. Chen,
$(K,K')$-quasiconformal harmonic mappings of the upper half plane onto itself,
\textit{Ann. Acad. Sci. Fenn. Math.}, {\bf 37} (2012), 265--276.

\bibitem{CPW} Sh. Chen, S. Ponnusamy and X. Wang,
Bloch constant and Landau's theorem for planar $p$-harmonic mappings,
\textit{J. Math. Anal. Appl.}, {\bf 373} (2011), 102--110.

\bibitem{CV} Sh. Chen and M. Vuorinen,
Some properties of a class of elliptic partial differential operators,
\textit{J. Math. Anal. Appl.}, {\bf 431} (2015), 1124--1137.


\bibitem{go1} F. W. Gehring and B. G. Osgood,
Uniform domains and the quasihyperbolic metric.
\textit{J. Analyse Math.} {\bf 36} (1979), 50--74.

\bibitem{go2} F. W. Gehring and B. P. Palka,
Quasiconformally homogeneous domains.
\textit{J. Analyse Math.} {\bf 30} (1976), 172--199.


\bibitem{H1} L. H\"{o}rmander,
An introduction to complex analysis in several variables,
3rd ed., xii+254 pp. North-Holland Mathematical Library 7. Amsterdam: North-Holland Publishing Co.,
1990.


\bibitem{K1} D. Kalaj,
Quasiconformal and harmonic mappings between Jordan domains,
\textit{Math. Z.}, {\bf 260} (2008), 237--252.

\bibitem{K2} D. Kalaj,
On quasiconformal harmonic maps between surfaces,
\textit{Int. Math. Res. Notices}, {\bf 2} (2015), 355--380.



\bibitem{KM0} D. Kalaj and M. Mateljevi\'{c},
Inner estimate and quasiconformal harmonic maps between smooth domains,
\textit{J. Anal. Math.}, {\bf 100} (2006), 117--132.


\bibitem{KM1} D. Kalaj and M. Mateljevi\'{c},
On certain nonlinear elliptic PDE and quasiconfomal maps between Euclidean surfaces,
\textit{Potential Anal.}, {\bf 34} (2011), 13--22.

\bibitem{KP0} D. Kalaj and M. Pavlovi\'{c},
Boundary correspondence under harmonic quasiconformal  diffeomorphisms
of a half-plane,
\textit{Ann. Acad. Sci. Fenn. Math.}, {\bf 30} (2005), 159--165.

\bibitem{KP} D. Kalaj and M. Pavlovi\'{c},
On quasiconformal self-mappings of the unit disk satisfying Poisson equation,
\textit{Trans. Amer. Math. Soc.}, {\bf 363} (2011), 4043--4061.

\bibitem{MKM} M. Kne$\check{z}$evi\'{c} and M. Mateljevi\'{c},
On the quasi-isometries of harmonic quasiconformal mappings,
\textit{J. Math. Anal. Appl.}, {\bf 334} (2007), 404--413.

\bibitem{LCW} P. Li, J. Chen and X. Wang,
Quasiconformal solutions of Poisson equations,
\textit{Bull. Aust. Math. Soc.}, {\bf 92} (2015), 420--428.


\bibitem{M} O. Martio,
On harmonic quasiconformal mappings,
\textit{ Ann. Acad. Sci. Fenn. Math.}, {\bf 425} (1968), 3--10.

\bibitem{MV} M. Mateljevi\'{c} and M. Vuorinen,
On harmonic quasiconformal quasi-isometries,
\textit{J. Inequal. Appl.}, (2010). Art. ID 178732.

\bibitem{MC} J. Mu and X. Chen,
Landau-type theorems for solutions of a quasilinear differential equation,
\textit{J. Math. Study}, {\bf 47}(3) (2014), 295--304.

\bibitem{O} A. Olofsson,
Differential operators for a scale of Poisson type kernels in the unit disc,
\textit{J. Anal. Math.}, {\bf 27} (2002), 365--372.


\bibitem{P} M. Pavlovi\'{c},
Boundary correspondence under harmonic quasiconformal homeomorphisms
of the unit disk,
\textit{Ann. Acad. Sci. Fenn. Math.}, {\bf 27}(2) (2002), 365--372.

\bibitem{P1} M. Pavlovi\'{c},
Decompositions of $L^p$ and Hardy spaces of polyharmonic functions,
\textit{J. Math. Anal. Appl.}, {\bf 216} (1997), 499--509.

\bibitem{si} S. Simi\'{c},
Lipschitz continuity of the distace ratio metric on the unit disk.
{\it Filomat} {\bf 27}(8) (2013), 1505--1509.

\bibitem{sivw} S. Simi\'{c},  M. Vuorinen, and G. Wang,
Sharp Lipschitz constants for the distance ratio metric.
\textit{Math. Scand.}  {\bf 116}(1) (2015), 89--103.


\bibitem{V} M. Vuorinen,
Exceptional sets and boundary behavior of quasiregular mappings in
$n$-space.
\textit{Ann. Acad. Sci. Fenn. Ser. A I Math.} Dissertationes {\bf 11} (1976), 1--44.

\bibitem{vo2} M. Vuorinen,
Conformal invariants and quasiregular mappings.
\textit{J. Analyse Math.} {\bf 45} (1985), 69--115.


\end{thebibliography}
\end{document}